\documentclass[12pt,a4paper]{article}
\usepackage{amsmath}
\usepackage{amssymb}
\usepackage{amsxtra}
\usepackage{latexsym}

\begin{document}
\count0=1


\centerline{\bf\Large Pointwise estimates of a solution}
\vskip10pt\centerline{\bf\Large of some boundary value problem}
\vskip10pt\centerline{\bf\Large for the Emden-Fowler equation}
\vskip20pt \centerline{\bf Burskii I.V.} \vskip10pt
\centerline{\small Donetsk National University, Ukraine}

\vskip10pt










\centerline{\large\bf Introduction} 
\vskip 0.5cm

In this paper we will consider the Dirichlet problem for the
well-known Emden-Fowler equation. It will obtain several pointwise
estimates of a solution of the problem:

\begin{equation}\label{bur:eq1}
-\Delta u+|u|^{q-1}u=0,\ \ q>1,\ \ x\in \Omega=B_1\backslash
B_d\end{equation}
\begin{equation}\label{bur:eq2}
u=0\ \ \ \ \ \ \ \ \ \ \ |x|=1\end{equation}
\begin{equation}\label{bur:eq3}u=k\ \ \ \ \ \ \ \ \ \ \ |x|=d\end{equation}
for the case $n>2,\ \ 1<q<\frac{n+2}{n-2}$.
\vskip10pt

The similar equation will be considered

\begin{equation}\label{bur:eq4}\Delta u+|u|^{q-1}u=0,\ \ q>1,\ \
x\in \Omega=B_1\backslash B_d\end{equation}
and the Laplace equation for a comparison

\begin{equation}\label{bur:eq5}\Delta u=0\end{equation}
with the same boundary value conditions (\ref{bur:eq2}) and
(\ref{bur:eq3}).

The first part of the paper is devoted to some general statements
on existence and uniqueness of the solution. In the second part we
obtain estimates of a solution that follow from the definition of
a generalization solution and classical inequalities.  By means of
the same arguments there was obtained a comparison of estimate for
solution of the problem $(\ref{bur:eq1}), (\ref{bur:eq2}),
(\ref{bur:eq3})$ with estimate for solution of the problem
$(\ref{bur:eq5}), (\ref{bur:eq2}), (\ref{bur:eq3})$. The main
theorem of this part is:

\vskip 0.3cm

{\bf Theorem 1.} For $n>2,\ \ 1<q<\frac{n+2}{n-2}$

1) The solution of the problem (\ref{bur:eq1}), (\ref{bur:eq2}),
(\ref{bur:eq3}) satisfies the estimate:

$$u(x)   \le   (c_1k
+c_2 k^{q-1} d^{\,2} )\left( \frac {d}{|x|} \right)^{n-2}$$

\noindent

2) the solution of the problem (\ref{bur:eq5}), (\ref{bur:eq2}),
(\ref{bur:eq3}) and positive spherically symmetric solution of the
problem (\ref{bur:eq4}), (\ref{bur:eq2}), (\ref{bur:eq3})
satisfies the estimate:
$$u(x)   \le   ck \left( \frac {d}{|x|} \right)^{n-2}.$$


In the third part it is suggested a method that based on
comparison theorems and permit to obtain both upper and lower
estimates which have some explicit constants.

\vskip 0.3cm

{\bf Theorem 2.} For the same $n$ and $q$

1) the solution of the problem $(\ref{bur:eq1}), (\ref{bur:eq2}),
(\ref{bur:eq3})$ has the following estimates
$$C_1 x^{2-n} + C_2 x^2 + C_3\le u(x)\le k\frac {x^{2-n}-1}{d^{2-n}-1}
$$

Here $C_1,C_2,C_3$ are some explicit constants (see formulae
(\ref{bur:eq24}), (\ref{bur:eq25})) in the case of small $x$ and
$d$ coefficient $C_1$ has the following principal term:
$$\frac{k-\frac{k^q}{n(n-2)}}{d^{2-n}-1}.$$
The left side may have only restricted application, when $k$ is
not a big constant.

2) Positive spherical solution of the problem (\ref{bur:eq4}),
(\ref{bur:eq2}), (\ref{bur:eq3}) has the following estimate
$$u(x)\ge k\frac {x^{2-n}-1}{d^{2-n}-1}.
$$

If besides this solution is bounded above by constant $k$ then the
estimate holds

$$C'_1 x^{2-n} + C'_2 x^2 + C'_3\ge u(x)\ge k\frac {x^{2-n}-1}{d^{2-n}-1}.
$$
with constants $C'_1,C'_2,C'_3$. Here for small $d$ the
coefficient $C'_1$ has the principal term:
$$\frac{k+\frac{k^q}{n(n-2)}}{d^{2-n}-1}.$$


The estimates that were obtained in theorems and those of similar
to them have important significance for studying nonlinear
elliptic boundary value problems in domain with finely granulated
boundary (see \cite{1}). Equations (\ref{bur:eq1}) and
(\ref{bur:eq4}) have also important physical applications. They
can be met in astrophysics in the form of Emden equation and in
atomic physics in the form of the Fermi-Thomas equation. Remark
that the question of nonuniqueness for equation (\ref{bur:eq4})
was studied by S.Pohozhaev in \cite{4}.

\vskip 30pt \centerline {\large\bf 1. Some general remarks.}
\vskip 30pt

Further we will consider a wide enough class of solutions
$C^2(\bar\Omega )$ , that is spherically symmetric if it won't be
underlined otherwise.

Remark, at first, that there exists some generalized solution $u$
(see the equality (\ref{bur:eq7}) for the equation
(\ref{bur:eq1})) from the Sobolev space $W^{1,2}_0(\Omega)$ of
homogeneous Dirichlet problem for the case of
$1<q<\frac{n+2}{n-2}$ as it flows out from theorem 4.1 in \cite{1}
(In this case for $n>2$ we have $W^{1,2}_0(\Omega)\subset
L_q(\Omega)$). For nonhomogeneous boundary data one may carry on
the same arguments for a remainder $u-\chi$.

Remark, secondly, that any generalized solution of homogeneous
Dirichlet problem is arbitrarily smooth function because it
assumes a raising of the smoothness by virtue of the ellipticity.
Indeed, if we carry over the term $u^{q-1}|u|$ into the right-part
side of the equation and we will suppose it belongs to any space
then the same $u$ as the solution of equation from left-part side
with known right-part side will a function from more smooth space.

Note, more, that our boundary value problem assumes only unique
solution for the equation (\ref{bur:eq1}) or (\ref{bur:eq5}) by
virtue of usual arguments with the maximum principle. Therefore
each solution of the problems (\ref{bur:eq1}), (\ref{bur:eq2}),
(\ref{bur:eq3}) or (\ref{bur:eq1}), (\ref{bur:eq2}),
(\ref{bur:eq5}) is spherically symmetric solution which satisfies
the equation
$$u''+\frac {n-1}{r}u'=u^q.$$
because respective boundary value problem (see (\ref{bur:eq17}))
has unique solution.

Let us, at last, give a proof that our solution of (\ref{bur:eq1})
is necessarily a positive function. Indeed, let the function $u$
be negative at some point. Then it has a minimum at some point
$x_0$. At this point we have $\nabla u(x_0)=0,\Delta u(x_0)>0$ but
the right-side part of the equation is negative. Therefore $u>0$.



\vskip20pt

\centerline{\large\bf 2.\ The obtaining of pointwise estimates of
Laplace}\vskip 10pt \centerline{\large\bf equation solution and
Emden-Fouler equation}\vskip 10pt \centerline{\large\bf solution
by means of the Mozer's inequality.} \vskip 20pt


Here we will prove the theorem 1. We consider the boundary value
problem (\ref{bur:eq2}), (\ref{bur:eq3}) for the equation
(\ref{bur:eq1}) and for convenience below we will compare
calculations for equations (\ref{bur:eq5}). And then we will
remark the case of the equation (\ref{bur:eq4}). \vskip10pt

{\bf 1).} First, we will obtain some  auxiliary estimates of
expressions
\[\int \limits_\Omega{\left|\frac {\partial u}{\partial
x}\right|}^2 \,dx\ \mbox{and}\ \int \limits_\Omega
\left({\left|\frac {\partial u}{\partial x}\right|}^2 + u^q\right)
\,dx.\] For this goal let us write integral identities for the
equations (\ref{bur:eq5}) and (\ref{bur:eq1}) respectively that we
will write as equalities
\begin{equation}\label{bur:eq6}
\sum \limits_{i=1}^n \int \limits_\Omega \frac {\partial
u}{\partial x_i} \frac {\partial \phi }{\partial x_i} \,dx=0,
\end{equation}
\begin{equation}\label{bur:eq7}\int \limits_\Omega \left\{ \frac {\partial u}{\partial x_i}
\frac {\partial \phi }{\partial x_i} + u^{q-1} \phi \right\}
\,dx=0
\end{equation}
that hold for all $\phi \in W_0^{1,2} (\Omega )$ Here we denote by
$\phi \in W_0^{1,2} (\Omega )$ the usual Sobolev space of
functions in $\Omega$, that are square integrable with their first
derivatives and that equal to zero on the boundary.

Let us take the function $\phi$ as the following
\[\phi(x) =u(x)-k\omega \left(\frac {|x|}{d}\right)\ \ \ \ \ \ \
2d<1\] where $\omega(s)\in C^\infty_0({\bf R}^+),\ {\bf R}^+
=[0,+\infty)$ is a function with properties $\omega (s)=1,
s\in[0,1]$; $\omega (s)=0, s\in[2,\infty)$; $0<\omega (s)<1,
s\in(1,2)$.

Remark, the function $\phi(x)$ is vanished for $|x|=1,\,|x|=d$ and
$$\int \limits_\Omega {\left |\frac {\partial u}{\partial x}\right |}^2 \,dx =
k\int \limits_\Omega \sum \limits_i \frac {\partial u}{\partial
x_i} \frac {\partial \omega }{\partial x_i} \left(  \frac{|x|}{d}
\right) \,dx$$

It is obvious that the estimate
$$\left | \frac {\partial \omega }{\partial x_i}
\left(  \frac {|x|}{d}  \right) \right |= \left |\omega ^\prime
\left(\frac {|x|}{d}\right) \frac {1}{d} \frac {x_i}{|x|}\right |
\leq \frac {C_1}{d}$$ .


Let us use together with well-known inequality $ab \leq \epsilon
a^2 + \frac {b^2}{\epsilon }$:

\[k\int \limits_\Omega \sum \limits_i \frac {\partial u}{\partial
x_i} \frac {\partial \omega }{\partial x_i} \left(  \frac {|x|}{d}
\right) \,dx \le kn \int \limits_{d\le |x| \le 2d} \left| \frac
{\partial u}{\partial x}\right| \frac {C_1}{d} \,dx \le\] \[
\le\epsilon \int \limits_\Omega {\left |\frac {\partial
u}{\partial x}\right |}^2 \,dx + \frac {1}{\epsilon } k^2 n^2
\frac {C_1^2}{d^2} \int \limits_{d\le |x| \le 2d}\,dx.\]

From here we have (for the equation (\ref{bur:eq6}))
\begin{equation}\label{bur:eq8}
\int \limits_\Omega {\left |\frac {\partial u}{\partial x}\right
|}^2 \,dx \le C_2 k^2 d^{n-2} \end{equation}

Let us now obtain an estimate for $\int\limits_\Omega ( {\left
|\frac {\partial u}{\partial x}\right |}^2 + u^q ) \,dx,$ where
$u$ is a solution of (\ref {bur:eq7}). To this end we substitute
the same function $\phi(x) =u(x)-k\omega \left(\frac
{|x|}{d}\right)$ into integral identity (\ref {bur:eq7}):
$$\int \limits_\Omega \left( {\left |\frac {\partial u}{\partial x}\right |}^2 +
u^q \right) \,dx = k\int \limits_\Omega \sum \limits_{i=1}^n \frac
{\partial u}{\partial x_i} \frac {\partial \omega }{\partial x_i}
\left(  \frac {|x|}{d}  \right) \,dx + k\int \limits_\Omega
u^{q-1} \omega \left( \frac {|x|}{d} \right) \,dx \le$$
$$\le \int
\limits_{d\le |x| \le 2d} \left| \frac {\partial u}{\partial
x}\right| kn \frac {C_1}{d} \,dx + \int \limits_\Omega u^{q-1} k
\omega \left( \frac {|x|}{d} \right) \,dx. $$ We used here
previously obtained estimate for $\left |\frac {\partial \omega
}{\partial x_i} \left( \frac {|x|}{d} \right) \right |$.

Then for the first term we use the inequality $ab \le \frac
{a^2}{2} + \frac {b^2}{2}$ and for the second one the same
inequality but with $\epsilon$:
$$\frac {\epsilon ab}{\epsilon } \le
\frac {\epsilon^p a^p}{p}+\frac {b^{\,r}}{\epsilon^r r},$$ where
we take $p=q/(q-1),r=q$ and the constant $\epsilon$ is chosen so
that $\epsilon^p/p=1/2$. We continue the inequality:
\[\int \limits_{d\le |x| \le 2d} \left| \frac {\partial
u}{\partial x}\right| kn \frac {C_1}{d} \,dx + \int \limits_\Omega
u^{q-1} k \omega \left( \frac {|x|}{d} \right) \,dx \le\]
$$\le\frac {1}{2} \int \limits_{d\le |x| \le 2d}
{\left| \frac {\partial u}{\partial x}\right|}^2 \,dx + \frac {k^2
n^2}{2} \int \limits_{d\le |x| \le 2d} \frac {C_1^2}{d^2} \,dx +
\frac {1}{2} \int \limits_\Omega u^q \,dx + C_2 \int \limits_{d\le
|x| \le 2d} k^q \omega^q \,dx.$$


hence we have finally
\begin{equation}\label{bur:eq9} \int \limits_\Omega \left( {\left
|\frac {\partial u}{\partial x}\right |}^2 + u^q \right) \,dx \le
C_1 k^2 d^{n-2} + C_2 k^q d^n,
\end{equation}
where $C_1,C_2$ are some new constants.

{\bf 2).} Now we get over some auxiliary estimate of expression
$\int \limits_{\left\{ u<\mu \right\} } {\left |\frac {\partial
u}{\partial x}\right |}^2 \,dx $ where $\mu:\ 0<\mu<k$. Take as
$\phi$ the function $\phi = \mbox{min}\{u(x),\mu \} - \frac {\mu
}{k} u(x)$ and substitute it at first into (\ref{bur:eq6}) and
then into (\ref{bur:eq7}).

It is easy to see that the function $\phi(x)$ vanish for
$|x|=1,\,|x|=d$ and

$$
\frac{\partial}{\partial x}\ \mbox{min}\{u(x),\mu\}=\left\{
\begin{array}{l}
{\partial u}/{\partial x}, \mbox{\small \ if\ }\{u<\mu\},\\ \\ \ \
\ 0,\mbox{\small \qquad if\ }\{u>\mu\}
\end{array}\right.
$$

\[\int \limits_{\left\{ u<\mu \right\} } {\left |\frac {\partial
u}{\partial x}\right |}^2 \,dx = \frac {\mu }{k} \int
\limits_\Omega {\left |\frac {\partial u}{\partial x}\right |}^2
\,dx. \] By using (\ref{bur:eq8}) we obtain the necessary estimate
(for (\ref{bur:eq6}))
\begin{equation}\label{bur:eq10}\int \limits_{\left\{ u<\mu \right\} }
{\left |\frac {\partial u}{\partial x}\right |}^2 \,dx \le C_2 \mu
k d^{\,n-2}. \end{equation}

On the same way we can obtain the following estimate by means of
the substitution $\phi$ into (\ref{bur:eq7}) and then using
(\ref{bur:eq9}):
$$\int \limits_{\left\{ u<\mu \right\} } {\left |\frac {\partial
u}{\partial x}\right |}^2 \,dx + \int \limits_{\left\{ u<\mu
\right\} } u^q \,dx + \mu \int \limits_{\left\{ u<\mu \right\} }
u^{q-1} \,dx =$$
$$=\frac {\mu}{k} \int \limits_\Omega \left( {\left |\frac {\partial u}{\partial x}\right |}^2 +
u^q \right) \,dx \le \frac {\mu}{k} (C_2 k^2 d^{n-2} + C_3 k^q
d^n)$$ Thus, we come to the estimate for the equation
(\ref{bur:eq7})
\begin{equation}\label{bur:eq11}\int \limits_{\left\{ u<\mu \right\} }
{\left |\frac {\partial u}{\partial x}\right |}^2 \,dx + \int
\limits_{\left\{ u<\mu \right\} } u^q \,dx + \mu \int
\limits_{\left\{ u<\mu \right\} } u^{q-1} \,dx \le \mu (C_2 k
d^{n-2} + C_3 k^{q-1} d^n)
\end{equation}


{\bf 3).} In order to obtain a final estimate of the maximum
module of solution we will use Mozer's method (\cite{1}, \cite{2},
\cite{3}). It is valid the inequality
\begin{equation}\label{bur:eq12}\max \limits_{\Gamma_\rho }
{|u(x)|}^2 \le \frac {C}{\rho^n} \int \limits_{\tilde \Gamma (\rho
)} u^2(x)
\,dx, \end{equation}
where $\Gamma_\rho$ is a spherical layer between spheres of
radiuses $\rho+\epsilon$ and $\rho-\epsilon$, $\tilde\Gamma(\rho)$
is a wider spherical layer which encloses the layer $\Gamma_\rho$.
We will need the following inequality of a type by Fridrichs or
Poincare
\begin{equation}\label{bur:eq13}\int \limits_{B(\rho )} {|u(x)|}^p \le C
\rho^p \int \limits_{B(R)} {\left |\frac {\partial u}{\partial
x}\right |}^2 \,dx,\ \ \ \ \ \mbox{\ where\ }\rho <R
\end{equation}
Let us get over the main part of the calculations \[\int
\limits_{\Gamma (\rho )} u^2(x) \,dx = \int \limits_{\Gamma (\rho
)} {\left[ min\{u(x),\mu \} \right]}^2 \le \int \limits_{B(\rho
+\epsilon )} min^2 \{u(x),\mu \} \,dx \le\]
\begin{equation}\label{bur:eq14}\le C\rho^2 \int \limits_{\left\{ u<\mu \right\} }
{\left |\frac {\partial u}{\partial x}\right |}^2 \,dx, \ \ \ \ \
\mu =\max \limits_{\Gamma_\rho } u(x).
\end{equation}

From (\ref{bur:eq12}), (\ref{bur:eq14}) and (\ref{bur:eq10}) we
have
$$\mu^2 \le C \frac {\rho^2}{\rho^n} \mu k d^{n-2}$$
$$\mu \le Ck {\left( \frac {d}{\rho} \right)}^{n-2}$$
$$\max \limits_{|x|=\rho } |u(x)| \le
Ck {\left( \frac {d}{\rho} \right)}^{n-2}.$$ Because $\rho$ can be
an arbitrary number between $d$ and 1, for the Laplace equation
(\ref{bur:eq6}) one may read the following pointwise estimate
\begin{equation}\label{bur:eq15}
u(x)\le max|u(x)|\le Ck {\left( \frac {d}{|x|} \right)}^{n-2}.
\end{equation}

Now, if one do the same for the equation (\ref{bur:eq7}), i.e. if
one uses (\ref{bur:eq12}), (\ref{bur:eq14}) and (\ref{bur:eq11}),
then he obtains \[\mu \le C_1 k {\left( \frac {d}{\rho}
\right)}^{n-2} + C_2 k^{q-1} d^2 {\left( \frac {d}{\rho}
\right)}^{n-2}\]


and finally
\begin{equation}\label{bur:eq16}u(x) \le \max|u(x)| \le
C_1 k {\left( \frac {d}{|x|} \right)}^{n-2} + C_2 k^{q-1} d^2
{\left( \frac {d}{|x|} \right)}^{n-2}
\end{equation}

Consider at last the case of the equation (\ref{bur:eq4}). For
this goal let us write integral identity for the equation
(\ref{bur:eq4}) \[\int \limits_\Omega \left\{ \frac {\partial
u}{\partial x_i} \frac {\partial \phi }{\partial x_i} - u^{q-1}
\phi \right\} \,dx=0, \] that holds for all $\phi \in W_0^{1,2}
(\Omega ).$ Here we do the same steps 1),2) and 3) as in the case
of the equation (\ref{bur:eq5}) without any problems because the
term with $-u^q$ is nonpositive. Therefore we obtain the same
estimate as for the equation (\ref{bur:eq5}).

The theorem 1 is proved.


\vskip 20pt

\centerline{\large\bf 3.\ Some way of obtaining of upper and lower
estimates} \vskip 10pt \centerline{\large\bf of a solution of the
Dirichlet boundary value problem} \vskip 10pt
\centerline{\large\bf for the Emden-Fouler equation.} \vskip 20pt

We will consider a spherically symmetric solution of the problem
(\ref{bur:eq1}), (\ref{bur:eq2}), (\ref{bur:eq3}) or problem
(\ref{bur:eq4}), (\ref{bur:eq2}), (\ref{bur:eq3}) belonging to the
space $C^2(\bar\Omega)$ under condition $1<q<(n+2)/(n-2)$.
(see section 2). Firstly, let remark that if $U(x)=u(|x|)$ is the
existing smooth solution of the problem (\ref{bur:eq1}),
(\ref{bur:eq2}), (\ref{bur:eq3}) then the function $u(r)$
satisfies the equation
\[u''+\frac {n-1}{r}u'=|u|^{q-1}u.\] The proof flows out from an
usual transfer to spherical variables. By means of solution
positiveness our problem can be written as
\begin{equation}\label{bur:eq17}\left\{\begin{array}{l}
u''+\frac {n-1}{r} u'=u^q \\ u(1)=0,\ \ \ u(d)=k.
\end{array}\right.\end{equation}

Now let us formulate the main auxiliary result of this
point.\vskip 10pt

{\bf Lemma 1.} Let the functions $\hat u_1 $, $\hat u_2 $ be
positive functions and the functions $ u_1 $, $ u_2 $ be solutions
of the following problems

\begin{equation}\label{bur:eq18}
\left\{
\begin{array}{l}
u_1''+\frac {n-1}{r}u_1'=\hat u_1^q ,\\ \\ u_1(1)=0,\ \ \
u_1(d)=k,
\end{array}\right.
\end{equation}
\begin{equation}\label{bur:eq19}
\left\{
\begin{array}{l}
u_2''+\frac {n-1}{r}u_2'=\hat u_2^q ,\\ \\ u_2(1)=0,\ \ \
u_2(d)=k.
\end{array}\right.
\end{equation}

Then $\hat u_1 \le \hat u_2$ implies $u_2 \le u_1$.

{\bf Proof.} Let us consider the remainder of equations
$(\ref{bur:eq19})-(\ref{bur:eq18})$. After subtraction we will
have
$$
(u_2-u_1)''+\frac {n-1}{r} (u_2-u_1)'=\hat u_2^q-\hat u_1^q \ge 0,
$$
$$
u_2(1)-u_1(1)=0,\ \ u_2(d)-u_1(d)=0. 
$$
Let $u=u_2-u_1\ge 0$ at some point $r$. Then it has a maximum at
some point $r_0$. At this point we have $u'(r_0)=0,u''(r_0)<0$ but
the right-side part of the equation is positive. Therefore $u\le
0$ that is $u_2\le u_1$. \vskip 5pt

Let us pass to description of the above-mentioned method.
\vskip10pt

{\bf Corollary 1. } Let $u$ be a solution of the problem
(\ref{bur:eq17}) and $u_0$ be any function that $u_0\le u$. Let
$u_1$ be a solution of the problem
\begin{equation}\label{bur:eq20}\left\{\begin{array}{l}
u_1''+\frac {n-1}{r} u_1'=u_0^q \\ u_1(1)=0,\ \ \ u_1(d)=k.
\end{array}\right.\end{equation}
Then $u\le u_1$.

For proof we apply the lemma to the problems (\ref{bur:eq17}) and
(\ref{bur:eq20}).

\vskip10pt

{\bf Corollary 2. } Let $u$ be a solution of the problem
(\ref{bur:eq17}) and $u_1$ be any function that $u_1\ge u$. Let
$u_2$ be a solution of the problem
\begin{equation}\label{bur:eq21}\left\{\begin{array}{l}
u_2''+\frac {n-1}{r} u_2'=u_1^q \\ u_2(1)=0,\ \ \ u_2(d)=k.
\end{array}\right.\end{equation}
Then $u\ge u_2$. \vskip10pt

For proof we apply the lemma 1 to the problems (\ref{bur:eq17})
and (\ref{bur:eq21}).

Let us consider our second problem with equation (\ref{bur:eq4})
that can be written as
$$\left\{\begin{array}{l} u''+\frac {n-1}{r} u'=-u^q \\ u(1)=0,\ \ \ u(d)=k.
\end{array}\right. \eqno(17^\prime)$$


{\bf Lemma 2.} Let the functions $\hat u_1 $, $\hat u_2 $ be
positive functions and the functions $ u_1 $, $ u_2 $ be solutions
of the following problems

$$
\left\{
\begin{array}{l}
u_1''+\frac {n-1}{r}u_1'=-\hat u_1^q ,\\ \\ u_1(1)=0,\ \ \
u_1(d)=k
\end{array}\right.$$
$$
\left\{
\begin{array}{l}
u_2''+\frac {n-1}{r}u_2'=-\hat u_2^q ,\\ \\ u_2(1)=0,\ \ \
u_2(d)=k
\end{array}\right.$$

Then $\hat u_1 \le \hat u_2$ implies $u_1 \le u_2$.

The proof is analogous to the proof of the lemma 1. \vskip10pt

{\bf Corollary 3. } Let $u$ be a solution of the problem
$(17^\prime)$ and $u_0$ be any function which $u_0\le u$. Let
$u_1$ be a solution of the problem
$$\left\{\begin{array}{l}
u_1''+\frac {n-1}{r} u_1'=-u_0^q, \\ u_1(1)=0,\ \ \ u_1(d)=k.
\end{array}\right.
$$
Then $u_1\le u$. \vskip10pt


{\bf Corollary 4. } Let $u$ be a solution of the problem
$(17^\prime)$ and $u_1$ be any function which $u_1\ge u$. Let
$u_2$ be a solution of the problem
$$\left\{\begin{array}{l}
u_2''+\frac {n-1}{r} u_2'=-u_1^q \\ u_2(1)=0,\ \ \ u_2(d)=k.
\end{array}\right.
$$
Then $u_2\ge u$. \vskip10pt

Thus, we can obtain a set of estimates of the solution $u$, lower
and upper.

For further applications we do some calculations. Consider the
problem
$$\left\{\begin{array}{l}
w''+\frac {n-1}{x} w'=f(x)\\ w(1)=0,\ \ \ w(d)=k.
\end{array}\right.$$
If denote $f(x)=u^q, w=u$ then equality (\ref{bur:eq17}) can be
written as such equation. Solving such ordinary linear
differential equation and boundary value problem as usually, we
receive
\begin{equation}\label{bur:eq22}w(x)=\int \limits_x^1 t^{1-n} \int \limits_t^1 \tau ^{n-1}
f(\tau ) \, d\tau dt  +  (k-D)\frac {x^{2-n}-1}{d^{2-n}-1},
\end{equation}
where \ \ $D=\int \limits_d^1
t^{1-n} \int \limits_t^1 \tau ^{n-1} 
f(\tau ) \, d\tau dt.$

\vskip10pt

Bring some calculations based on this approach.

\vskip10pt

I 1). Let $u_0=0$ for the equation (\ref{bur:eq1}). Because
$u_0\le u$ we have $u\le u_1$ from the corollary 1, that
\begin{equation}\label{bur:eq23}u(x)\le k\frac
{x^{2-n}-1}{d^{2-n}-1}\end{equation}

I 2) Now for the equation (\ref{bur:eq1}) we would like to take
$u_1=k$ because $u\le k$. We will receive

$$u_2 = \int \limits_x^1 t^{1-n} \int \limits_t^1 \tau ^{n-1} k^q \, d\tau dt +
(k-D_0)\frac {x^{2-n}-1}{d^{2-n}-1},$$
$$D_0 = \int \limits_d^1 t^{1-n} \int \limits_t^1 \tau ^{n-1} k^q \, d\tau dt =
\frac {k^q}{n} \left[ \frac {1}{2-n} - \frac {d^{2-n}}{2-n} -
\frac {1}{2} + \frac {d^2}{2} \right]=$$
$$= \frac {k^q}{2n(n-2)} (2d^{2-n}+ (n-2)d^2 -n).$$


Denote \begin{equation}\label{bur:eq24}\qquad\qquad\tilde C_1 =
\frac {k^q}{2n(n-2)},\ \ \ \tilde C_2 = -\frac {k-D_0}{d^{2-n}-1},
\end{equation}

then
$$
u_2 = \tilde C_1 (2x^{2-n} + (n-2)x^2 - n) + \tilde C_2
(1-x^{2-n})= \qquad\qquad\qquad\ \ \ $$
\begin{equation}\label{bur:eq25}\qquad\qquad
(2\tilde C_1 - \tilde C_2) x^{2-n} +  \tilde C_1
(n-2)x^2 +(\tilde C_2 - n\tilde C_1) =
 C_1 x^{2-n} + C_2 x^2 +
C_3\end{equation}

and corollary 2 says
$$u \ge C_1 x^{2-n} + C_2 x^2 + C_3. $$
Here for small $d$ the coefficient $C_1=2\tilde C_1-\tilde C_2$
has the following principal term:
$$\frac{k-\frac{k^q}{n(n-2)}}{d^{2-n}-1},$$
therefore this formula can have only bounded application, i.e.
when $k$ is not a big constant. \vskip 20pt

II 1) Let $u_0=0$ for equation (\ref{bur:eq4}). Because it is
considered a positive solution i.e. $u_0\le u$ we have $u_1\le u$
from the corollary 3 or
$$u(x)\ge k\frac {x^{2-n}-1}{d^{2-n}-1}
$$

II 2) Now for the equation (\ref{bur:eq4}) we assume that $u\le k$
and take $u_1=k$. We will receive the same formula as in the point
I 2) but with the change of $k^q$ on $-k^q$. Then corollary 4 says
$$u \le C'_1 x^{2-n} + C'_2 x^2 + C'_3 $$
and for small $d$ the coefficient $C'_1$ has the principal term:
$$\frac{k+\frac{k^q}{n(n-2)}}{d^{2-n}-1}.$$

The proof of the theorem 2 is over.

\end{document}